\documentclass[12pt]{article}
\usepackage{amssymb}
\usepackage{latexsym,bm}
\usepackage{graphicx}
\usepackage{amsmath}
\usepackage{mathrsfs}

\date{}
\begin{document}
\title{On the metric subgraphs of a graph\footnote{E-mail addresses:
{\tt huyanan530@163.com}(Y.Hu),
{\tt zhan@math.ecnu.edu.cn}(X.Zhan).}}
\author{\hskip -10mm Yanan Hu and Xingzhi Zhan\thanks{Corresponding author.}\\
{\hskip -10mm \small Department of Mathematics, East China Normal University, Shanghai 200241, China}}\maketitle
\begin{abstract}
The three subgraphs of a connected graph induced by the center, annulus and periphery are called its metric subgraphs. The main results are as follows.
(1) There exists a graph of order $n$ whose metric subgraphs are all paths if and only if $n\ge 13$ and the smallest size of such a graph of order $13$ is $22;$ (2) there exists a graph of order $n$ whose metric subgraphs are all cycles if and only if $n\ge 15,$ and there are exactly three such graphs of order $15;$ (3) for every integer $k\ge 3,$ we determine the possible orders for the existence of a graph whose metric subgraphs are all connected $k$-regular graphs; (4) there exists a graph of order $n$ whose metric subgraphs are connected and pairwise isomorphic if and only if $n\ge 24$ and $n$ is divisible by $3.$ An unsolved problem is posed.
\end{abstract}



\section{Introduction}

We consider finite simple graphs. For terminology and notations we follow the books [3] and [11]. The {\it order} of a graph $G$, denoted $|G|,$ is its number
of vertices, and the {\it size} is its number of edges. We denote by $V(G)$ and $E(G)$ the vertex set and edge set of a graph $G$ respectively. Denote by $d_{G}(u,v)$ the distance between two vertices $u$ and $v$ in $G.$ If the graph $G$ is clear from the context, we simply write $d(u,v).$ The {\it eccentricity}, denoted by $e(v),$ of a vertex $v$ in a graph $G$ is the distance to a vertex farthest from $v.$ Thus $e(v)={\rm max}\{d(v,u)| u\in V(G)\}.$ If $e(v)=d(v,x),$ then the vertex $x$ is called an {\it eccentric vertex} of $v.$ The {\it radius} of a graph $G,$ denoted ${\rm rad}(G),$ is the minimum eccentricity of all the vertices in $V(G),$ whereas the {\it diameter} of $G,$ denoted ${\rm diam}(G),$ is the maximum eccentricity. For graphs we will use equality up to isomorphism, so $G_1=G_2$ means that $G_1$ and $G_2$ are isomorphic. A graph is called {\it null} if it has order $0;$ otherwise it is {\it non-null.}

Let $G$ be a connected graph. A vertex $u$ is a {\it central vertex} of $G$ if $e(u)={\rm rad}(G).$ The {\it center} of $G,$ denoted $C(G),$ is the set of all central vertices of $G.$ A vertex $v$ is a {\it peripheral vertex} of $G$ if $e(v)={\rm diam}(G).$ The {\it periphery} of $G,$ denoted $P(G),$ is the set
of all peripheral vertices. A vertex $w$ is an {\it annular vertex} of $G$ if ${\rm rad}(G)<e(w)< {\rm diam}(G).$ The {\it annulus} of $G,$ denoted $A(G),$ is the set of all annular vertices.

{\bf Definition 1.} Let $G$ be a connected graph. The subgraph of $G$ induced by its center is called the {\it central subgraph} of $G;$
the subgraph of $G$ induced by its annulus is called the {\it annular subgraph} of $G;$ the subgraph of $G$ induced by its periphery is called the
{\it peripheral subgraph} of $G.$ These three subgraphs are called the {\it metric subgraphs} of $G.$

A self-centered graph has empty annulus, and hence its annular graph is the null graph. There are much work studying the central subgraph (e.g. [2],
[3, Chapter 2], [7], [8]),  some studying the peripheral subgraph ([1], [3], [4])  and little studying the annular subgraph [5]. In this paper, we consider these three subgraphs as a whole.

We will determine possible orders for the existence of a graph whose metric subgraphs are either all paths, or all cycles, or all connected $k$-regular graphs with $k\ge 3,$ or connected and pairwise isomorphic graphs,  and we also consider the smallest size problem. At the end we pose an unsolved
problem.

\section{Main results}

For two graphs $G$ and $H,$ $G\vee H$ denotes the {\it join} of $G$ and $H,$ which is obtained from the disjoint union $G+H$
by adding edges joining every vertex of $G$ to every vertex of $H.$  Given two vertex subsets $S$ and $T$ of a graph, we denote
by $[S, \,T]$ the set of edges having one endpoint in $S$ and the other in $T.$ $P_n,$ $C_n$ and $K_n$ denote the path of order $n,$ the cycle of order $n$
and the complete graph of order $n$ respectively. As usual, $q K_2$ denotes the graph consisting of $q$ pairwise vertex-disjoint edges, and ${\rm deg}(v)$ denotes the degree of a vertex $v.$ We denote by $\overline{G}$ the complement of a graph $G.$

{\bf Definition 2.} Given a sequence of graphs $H_1, H_2,\ldots, H_p$, their {\it circular join} is defined to be the graph obtained from the disjoint union
$H_1+H_2+\ldots+H_p$ by adding edges joining each vertex of $H_i$ to each vertex of $H_{i+1}$ for every $i=1,2,\ldots,p$ where $H_{p+1}$ means $H_1.$

{\bf Definition 3.} Let $G$ and $H$ be two graphs with $|G|\le |H|.$ Labeling the vertices of $G$ and $H$ by $x_1,\ldots, x_s$ and $y_1,\ldots, y_t$
respectively, the {\it nice connection} of $G$ and $H$ with respect to this vertex labeling is the graph obtained from the disjoint union $G+H$
by adding the edges $x_iy_i,$ $i=1,2,\ldots, s$ and the edges $x_1y_j,$ $j=s+1,\ldots,t.$

Given two graphs, there are possibly many nice connections of them, depending on the vertex labeling. For our purposes below, any nice connection
works.

We will need the following lemmas.

{\bf Lemma 1} (Lesniak [10]). {\it Let $G$ be a connected graph of order $n.$ Then for every integer $k$ with ${\rm rad}(G)<k\le {\rm diam}(G),$ there exist
at least two vertices in $G$ of eccentricity $k.$}

{\bf Lemma 2.} {\it If $G$ is a graph with a nonempty annulus, then ${\rm rad}(G)\ge 2,$ ${\rm diam}(G)\ge 4,$ and $|A(G)|\ge 2.$}

{\bf Proof.} Since ${\rm diam}(G)\le 2\, {\rm rad}(G)$ [11, p.78], if ${\rm rad}(G)=1$ then ${\rm diam}(G)\le 2,$ implying that the annulus of $G$ is empty,
a contradiction. Hence ${\rm rad}(G)\ge 2.$ Consequently ${\rm diam}(G)\ge {\rm rad}(G)+2\ge 4.$

The assertion $|A(G)|\ge 2$ follows from Lemma 1 and the condition that $G$ has a nonempty annulus. $\Box$

We denote by $e_G(v)$ the eccentricity of a vertex $v$ in $G.$ Recall that $|G|$ denotes the order of a graph $G.$ The following lemma is of independent
interest.

{\bf Lemma 3.} {\it Let $H$ be the peripheral subgraph of a connected graph $G.$ If $H$ is connected, then ${\rm rad}(H)\ge {\rm diam}(G)$
and $|H|\ge 2\, {\rm diam}(G).$}

{\bf Proof.} Let $v$ be a central vertex of $H$ and let $x$ be an eccentric vertex of $v$ in $G.$ Then $x\in P(G)=V(H).$ We have
$$
{\rm diam}(G)=e_G(v)=d_G(v,x)\le d_H(v,x)\le e_H(v)={\rm rad}(H),
$$
showing that ${\rm rad}(H)\ge {\rm diam}(G).$ Combining this inequality with the fact that ${\rm rad}(H)\le |H|/2$ we obtain $|H|\ge 2\, {\rm diam}(G).$ $\Box$

For a graph $G$ and $S\subseteq V(G),$ the {\it neighborhood} of $S$ is defined to be $N(S)=\{x\in V(G)\setminus S|\, x\,\, {\rm has}\,\,{\rm a}\,\,{\rm neighbor}\,\,
 {\rm in}\,\, $S$\}.$  We will repeatedly use the fact that the eccentricities of two adjacent vertices differ by at most $1.$

{\bf Lemma 4.} {\it Let $W$ be the annular subgraph of a connected graph $G.$ If $W$ is non-null and connected, then ${\rm rad} (W)\ge 2$ and
consequently $|A(G)|\ge 4.$ }

{\bf Proof.} Clearly $N(C(G))\subseteq A(G).$ To the contrary, suppose ${\rm rad} (W)=1.$ Let $v$ be a central vertex of $W.$
Denote $r={\rm rad}(G).$ Let $x$ be any vertex in $V(G).$ If $x\in C(G)$ then $d(v,x)\le r;$ if $x\in A(G)$ then $d(v,x)\le 1\le r.$
If $x\in P(G),$ choose any vertex $y\in C(G)$ and let $Q=y,...,z,...,x$ be a shortest $(y,x)$-path in $G$ where $z\in A(G).$
Then $Q$ has length at most $r$ and the subpath $Q[z, x]$ has length at most $r-1.$ Again we have $d(v,x)\le d(v,z)+d(z,x)\le 1+(r-1)\le r.$ This shows that $e(v)\le r,$ a contradiction. Hence ${\rm rad} (W)\ge 2.$

Since a connected graph of order at most $3$ has radius at most $1,$ we obtain $|A(G)|\ge 4.$ $\Box$

{\bf Lemma 5.} {\it If a connected graph has a connected peripheral subgraph and a non-null connected annular subgraph, then its order is at least $13.$ }

{\bf Proof.} Suppose $G$ is a connected graph whose peripheral subgraph is connected and whose annular subgraph is  non-null and connected.
Consider the diameter. By Lemma 2, ${\rm diam}(G)\ge 4.$ Then by Lemma 3, $|P(G)|\ge 2\,{\rm diam}(G)\ge 8.$ Lemma 4 gives $|A(G)|\ge 4.$ Note that every non-null graph has at least one central vertex. We obtain
$$
|G|=|C(G)|+|A(G)|+|P(G)|\ge 1+4+8=13.
$$
This completes the proof. $\Box$

{\bf Theorem 6.} {\it There exists a connected graph of order $n$ whose metric subgraphs are all paths if and only if $n\ge 13,$ and the smallest size of such a graph of order $13$ is $22.$}

{\bf Proof.} If $G$ is a connected graph of order $n$ whose metric subgraphs are all paths, then by Lemma 5, $n\ge 13.$

Conversely, for every $n\ge 13$ we will construct such a graph of order $n.$ First, the graph $G_1$ in Figure 1 is a graph of order $13$ whose metric subgraphs are all paths.
\vskip 3mm
\par
 \centerline{\includegraphics[width=2.6in]{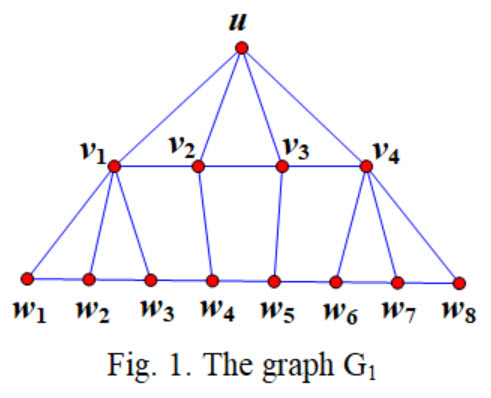}}
\par
The central subgraph, the annular subgraph and the peripheral subgraph of $G_1$ are the paths $u,$ $Q=v_1v_2v_3v_4,$ and $w_1w_2...w_8$ respectively. Next for a given integer $n\ge 14,$ in $G_1$ replacing the vertex $u$ by the path $P_{n-12}$ and then taking the join of $P_{n-12}$ and the path $Q$ we obtain a graph.
This is a graph of order $n$ whose metric subgraphs are all paths.

Finally we show that the smallest size of such a graph of order $13$ is $22.$ Let $H$ be a graph of order $13$ whose metric subgraphs are all paths.
Denote the center, annulus and periphery of $H$ by $C,$ $A$ and $P$ respectively. By the proof of Lemma 5, we have
$$
|C|=1,\,\,\, |A|=4,\,\,\, |P|=8,\,\,\,{\rm diam}(H)=4\,\,\,{\rm and}\,\,\,{\rm hence}\,\,\, {\rm rad}(H)=2.
$$
Let $C=\{u\},$ let the annular subgraph of $H$ be the path $v_1v_2v_3v_4,$ and let the peripheral subgraph of $H$ be the path $w_1w_2...w_8.$
Then $e(u)=2,$ $e(v_i)=3,\, i=1,\ldots, 4,$ and $e(w_j)=4,\, j=1,\ldots,8.$ Since the eccentricities of two adjacent vertices differ by at most $1,$
$N(u)\subseteq A.$ Now the condition $e(u)=2$ implies that each vertex in $P$ has a neighbor in $A.$ Hence $|[P,\, A]|\ge 8.$ Clearly, any eccentric vertex of a vertex in $P$ lies in $P.$ Since $e(v_2)=3,$ every eccentric vertex of $v_2$ lies in $P.$ Let $w_k$ be an eccentric vertex of $v_2.$ Then $v_4$ is the unique
neighbor of $w_k$ in $A.$ Since $d(u,w_k)\le 2,$ $u$ and $v_4$ are adjacent. Similarly, considering the vertex $v_3$ we deduce that $u$ and $v_1$ are adjacent.
We claim that $N(u)=A.$ Otherwise ${\rm deg}(u)\le 3.$ Without loss of generality, suppose $v_3$ and $u$ are nonadjacent. Denote $S=\{v_1,\, v_2\}$ and
$T=\{v_4\}.$ Then each vertex in $P$ has a neighbor in $S\cup T.$ We will repeatedly use this fact. Since $w_1$ is the only possible eccentric vertex of $w_5,$
we have $d(w_1, w_5)=4.$ Let $x\in S\cup T$ be a neighbor of $w_1$ and let $y\in S\cup T$ be a neighbor of $w_5.$ Then we have the following two possible cases.

Case 1. $x\in S$ and $y\in T.$ Since $w_8$ is the unique eccentric vertex of $w_4,$ $w_8$ and $v_4$ are nonadjacent. Hence $w_8$ has a neighbor in $S.$
Consequently, the neighbor of $w_4$ in $S\cup T$ must be $v_4.$ To keep $d(w_4, w_8)=4,$ $w_7$ and $v_4$ cannot be adjacent. Thus $w_7$ has a neighbor in
$S.$ Similarly, to keep $d(w_5, w_1)=4,$ $w_2$ and $v_4$ cannot be adjacent. Thus $w_2$ has a neighbor in $S.$ Now $w_6$ is the only eccentric vertex
of $w_2.$ Hence $w_6$ has no neighbor in $S,$ implying that $v_4$ is the only neighbor of $w_6$ in $S\cup T.$ Since $w_3$ has a neighbor in $S\cup T,$
we deduce that $e(w_7)\le 3,$ a contradiction.

Case 2. $x\in T$ and $y\in S.$ The condition $e(w_5)=4$ implies that the two vertices $w_4$ and $w_6$ are nonadjacent to $v_4.$ Thus, both $w_4$ and $w_6$
have a neighbor in $S.$ Since $w_8$ is the unique eccentric vertex of $w_4,$ $w_8$ and $v_4$ are adjacent and $w_3$ has a neighbor in $S.$ If $w_7$
is adjacent to $v_4,$ using the fact that $w_2$ has a neighbor in $S\cup T$ we deduce that $e(w_6)\le 3,$ a contradiction. Hence $w_7$ has a neighbor
in $S.$ But then $e(w_7)\le 3,$ a contradiction again. This shows that ${\rm deg}(u)=4.$

We conclude that the size of $H$ is at least $3+7+8+4=22.$ Conversely, the graph $G_1$ in Figure 1 is a graph of order $13$ and size $22$ whose metric
subgraphs are all paths. This completes the proof. $\Box$

{\bf Remark 1.} By the proof of Theorem 6, it is not difficult to check that there are exactly $64$ connected graphs of order $13$ and size $22$
whose metric subgraphs are all paths.

We will need the following two results.

{\bf Lemma 7} (Kim, Rho, Song and Hwang [9]) {\it Let $G$ be a graph of order $n$ with radius $r$ and minimum degree $k$ where $r\ge 3$ and $k\ge 2.$
Then $n\ge 2r(k+1)/3.$}

{\bf Lemma 8.} {\it Let $k$ and $n$ be integers with $1\le k\le n-1.$ Then there exists a $k$-regular graph of order $n$ if and
only if $kn$ is even. If $kn$ is even and $k\ge 2,$ then there exists a hamiltonian $k$-regular graph of order $n.$ }

Lemma 8 can be found in [6,  pp.12-13]. Its first part is well-known, but its hamiltonian part is usually not stated.

Next we determine the possible orders for the existence of a graph whose metric subgraphs are all connected and $k$-regular.
The answer depends on the nature of $k$ and there are six cases.

{\bf Theorem 9.} {\it Let $k\ge 2$ be an integer and denote $q=\lfloor k/3\rfloor.$ There exists a connected graph of order $n$ whose metric subgraphs are all connected and $k$-regular if and only if (1) $n\ge 14q+6$ when $k\equiv 0\,\, {\rm mod}\,\,3$ and $q$ is even; (2) $n$ is even and $n\ge 14q+8$ when
$k\equiv 0\,\, {\rm mod}\,\,3$ and $q$ is odd; (3) $n$ is even and $n\ge 14q+12$ when $k\equiv 1\,\, {\rm mod}\,\,3$ and $q$ is even; (4) $n\ge 14q+11$ when $k\equiv 1\,\, {\rm mod}\,\,3$ and $q$ is odd; (5) $n\ge 14q+15$ when $k\equiv 2\,\, {\rm mod}\,\,3$ and $q$ is even; (6) $n$ is even and $n\ge 14q+16$ when $k\equiv 2\,\, {\rm mod}\,\,3$ and $q$ is odd.}

{\bf Proof.} Let $G$ be a connected graph of order $n$ whose metric subgraphs are all connected and $k$-regular. Denote the center, annulus and periphery of $G$
by $C,$ $A$ and $P$ respectively.

Since the central subgraph $G[C]$ is $k$-regular, we have $|C|\ge k+1.$

We then estimate the cardinality of $A.$  Since the annular subgraph $W=G[A]$ is $k$-regular, we have $|A|\ge k+1$ where equality holds if and only if
$W$ is complete. By Lemma 4, ${\rm rad}(W)\ge 2,$ implying that $W$ is incomplete. Hence $|A|\ge k+2.$  When $k$ is odd, the order of $W$ must be even
since $W$ is $k$-regular, and consequently $|A|\ge k+3.$ It follows that (1) $|A|\ge 3q+2$ when $k\equiv 0\,\, {\rm mod}\,\,3$ and $q$ is even;
(2) $|A|\ge 3q+3$ when $k\equiv 0\,\, {\rm mod}\,\,3$ and $q$ is odd; (3) $|A|\ge 3q+4$ when $k\equiv 1\,\, {\rm mod}\,\,3$ and $q$ is even;
(4) $|A|\ge 3q+3$ when $k\equiv 1\,\, {\rm mod}\,\,3$ and $q$ is odd; (5) $|A|\ge 3q+4$ when $k\equiv 2\,\, {\rm mod}\,\,3$ and $q$ is even;
(6) $|A|\ge 3q+5$ when $k\equiv 2\,\, {\rm mod}\,\,3$ and $q$ is odd.

Next we estimate the cardinality of $P.$ Denote by $H$ the peripheral subgraph of $G.$ By Lemma 2, ${\rm diam} (G)\ge 4$ and by Lemma 3,
${\rm rad}(H)\ge {\rm diam}(G).$ Thus $r={\rm rad}(H)\ge 4.$ By Lemma 7 we obtain $|P|\ge 2r(k+1)/3,$ from which we deduce the following information.
(i) $|P|\ge 8q+3$ when $k\equiv 0\,\, {\rm mod}\,\,3$ and $q$ is even; (ii)  $|P|\ge 8q+4$ when $k\equiv 0\,\, {\rm mod}\,\,3$ and $q$ is odd;
(iii) $|P|\ge 8q+6$ when $k\equiv 1\,\, {\rm mod}\,\,3;$ (iv)  $|P|\ge 8q+8$ when $k\equiv 2\,\, {\rm mod}\,\,3.$

Using the above estimation we can obtain a lower bound for $n=|C|+|A|+|P|.$ We have (1) $n\ge 14q+6$ when $k\equiv 0\,\, {\rm mod}\,\,3$ and $q$ is even; (2) $n$ is even and $n\ge 14q+8$ when $k\equiv 0\,\, {\rm mod}\,\,3$ and $q$ is odd; (3) $n$ is even and $n\ge 14q+12$ when $k\equiv 1\,\, {\rm mod}\,\,3$ and $q$ is even; (4) $n\ge 14q+11$ when $k\equiv 1\,\, {\rm mod}\,\,3$ and $q$ is odd; (5) $n\ge 14q+15$ when $k\equiv 2\,\, {\rm mod}\,\,3$ and $q$ is even; (6) $n$ is even and $n\ge 14q+16$ when $k\equiv 2\,\, {\rm mod}\,\,3$ and $q$ is odd.

Conversely, for every order $n$ in the range as stated in the theorem, we construct a connected graph of order $n$ whose metric subgraphs are all connected and $k$-regular.

Notation. For a positive even integer $f,$ we denote by $K_f-PM$ the graph obtained from $K_f$ by deleting a perfect matching; i.e., $K_f-PM=\overline{(f/2)K_2}.$ For an integer $g\ge 3,$ we denote by $K_g-HC$ the graph $\overline{C_g}.$ For an integer $s\ge 2,$ we denote by
$K_s-HP$ the graph $\overline{P_s}.$

Case (1). $k\equiv 0\,\, {\rm mod}\,\,3$ and $q$ is even. Now $k=3q$ is even. Let $n\ge 14q+6.$ We have $n-(11q+5)\ge k+1.$ By Lemma 8, there exists a connected
$k$-regular graph $C^{(1)}$ of order $n-(11q+5).$  We denote by $P^{(1)}$ the circular join of the eight graphs $H_1^{(1)},\ldots, H_8^{(1)}$ where
$H_1^{(1)}=H_4^{(1)}=H_7^{(1)}=K_{q+1},$ $H_2^{(1)}=H_3^{(1)}=H_5^{(1)}=H_6^{(1)}=K_q$ and $H_8^{(1)}=K_q-PM.$ Denote $A^{(1)}=K_{3q+2}-PM.$ Partition
the vertex set of $A^{(1)}$ into eight subsets $V_1^{(1)},\ldots, V_8^{(1)}$ such that $A^{(1)}[V^{(1)}_i]=K_1$ for $i=1,5,$ $A^{(1)}[V^{(1)}_j]=K_{q/2}$
for $j=2,3,4,6,7,8$ and  $A^{(1)}[V^{(1)}_1\cup V^{(1)}_5]=\overline{K_2},$ $A^{(1)}[V^{(1)}_s\cup V^{(1)}_{s+4}]=K_q-PM$ for $s=2,3,4.$ Let
$R^{(1)}_i=A^{(1)}[V^{(1)}_i],$ $i=1,\ldots, 8.$ Finally let $M_1$ be the graph obtained from the disjoint union $C^{(1)}+A^{(1)}+P^{(1)}$ by first
taking a nice connection of $R^{(1)}_i$ and $H^{(1)}_i$ for $i=1,\ldots, 8$ and then adding edges joining every vertex of $C^{(1)}$ to every vertex of
$A^{(1)}.$

Note that every vertex in $H^{(1)}_i$ has a unique neighbor in $R^{(1)}_i$ for $i=1,\ldots,8.$ It is easy to verify that ${\rm rad}(M_1)=2$ and
${\rm diam}(M_1)=4,$ the three graphs $C^{(1)},$ $A^{(1)}$ and $P^{(1)}$ are connected and $k$-regular, and they are the central subgraph, annular subgraph and
peripheral subgraph of $M_1$ which has order $n.$

Case (2). $k\equiv 0\,\, {\rm mod}\,\,3$ and $q$ is odd. Now $k=3q$ is odd. Let $n\ge 14q+8$ and $n$ is even.
We have $n-(11q+7)\ge k+1.$ By Lemma 8, there exists a connected $k$-regular graph $C^{(2)}$ of order $n-(11q+7).$
We denote by $P^{(2)}$ the circular join of the eight graphs $H_1^{(2)},\ldots, H_8^{(2)}$ where
$H_1^{(2)}=H_2^{(2)}=H_5^{(2)}=H_6^{(2)}=K_q,$ $H_3^{(2)}=H_4^{(2)}=H_7^{(2)}=H_8^{(2)}=K_{q+1}-PM.$
Denote $A^{(2)}=K_{3q+3}-HC.$ We distinguish two subcases.

Subcase (2.1). $q=1.$  In this case $A^{(2)}=K_6-HC.$ Let $V(A^{(2)})=\{u_1,u_2,u_3,u_4,u_5,u_6\}$ and
$E(A^{(2)})=\{u_iu_j|\, j\ne i+1\}$ where $u_7=u_1.$ Let $M_2$ be the graph obtained from the disjoint union $C^{(2)}+A^{(2)}+P^{(2)}$ by first
adding edges joining $u_i$ to each vertex in $H_j^{(2)}$ for $(i,j)=(1,1),\,(2,2),\,(2,3),\,(3,6),\,(3,7),\,(4,8),\,(5,4),\,(6,5),$
then adding edges joining every vertex of $C^{(2)}$ to every vertex of $A^{(2)}.$ It is easy to verify that ${\rm rad}(M_2)=2$ and
${\rm diam}(M_2)=4,$ the three graphs $C^{(2)},$ $A^{(2)}$ and $P^{(2)}$ are connected and $k$-regular, and they are the central subgraph, annular subgraph and
peripheral subgraph of $M_2$ which has order $n.$

Subcase (2.2). $q\ge 3.$ Let $V(A^{(2)})=\{u_1,u_2,\ldots, u_{3q+3}\}$ and $E(A^{(2)})=\{u_iu_j|\, j\ne i+1\}$ where $u_{3q+4}=u_1.$
Partition the vertex set of $A^{(2)}$ into eight subsets (four pairs)
\begin{align*}
&V_1^{(2)}=\{u_1\},\quad V_5^{(2)}=\{u_2\}\\
&V_2^{(2)}=\{u_{3+2j}|\,j=0,1,\ldots,(q-3)/2\},\quad V_6^{(2)}=\{u_{4+2j}|\,j=0,1,\ldots,(q-3)/2\}\\
&V_3^{(2)}=\{u_{q+2j}|\,j=1,\ldots,(q+1)/2\},\quad V_7^{(2)}=\{u_{q+1+2j}|\,j=1,\ldots,(q+1)/2\}\\
&V_4^{(2)}=\{u_{2q+1+2j}|\,j=1,\ldots,(q+1)/2\},\quad V_8^{(2)}=\{u_{2q+2+2j}|\,j=1,\ldots,(q+1)/2\}
\end{align*}
such that $A^{(2)}[V^{(2)}_j]=K_1$ for $j=1,5,$ $A^{(2)}[V^{(2)}_j]=K_{(q-1)/2}$
for $j=2,6,$ $A^{(2)}[V^{(2)}_j]=K_{(q+1)/2}$ for $j=3,4,7,8$ and  $A^{(2)}[V^{(2)}_1\cup V^{(2)}_5]=\overline{K_2},$
$A^{(2)}[V^{(2)}_2\cup V^{(2)}_{6}]=K_{q-1}-HP,$ $A^{(2)}[V^{(2)}_s\cup V^{(2)}_{s+4}]=K_{q+1}-HP$ for $s=3,4.$ Let
$R^{(2)}_i=A^{(2)}[V^{(2)}_i],$ $i=1,\ldots, 8.$ Finally let $M_2$ be the graph obtained from the disjoint union $C^{(2)}+A^{(2)}+P^{(2)}$ by first
taking a nice connection of $R^{(2)}_i$ and $H^{(2)}_i$ for $i=1,\ldots, 8$ and then adding edges joining every vertex of $C^{(2)}$ to every vertex of
$A^{(2)}.$

Note that every vertex in $H^{(2)}_i$ has a unique neighbor in $R^{(2)}_i$ for $i=1,\ldots,8.$ It is easy to verify that ${\rm rad}(M_2)=2$ and
${\rm diam}(M_2)=4,$ the three graphs $C^{(2)},$ $A^{(2)}$ and $P^{(2)}$ are connected and $k$-regular, and they are the central subgraph, annular subgraph and
peripheral subgraph of $M_2$ which has order $n.$

Case (3). $k\equiv 1\,\, {\rm mod}\,\,3$ and $q$ is even. Now $k=3q+1$ is odd. Let $n$ be even and $n\ge 14q+12.$ We have $n-(11q+10)\ge k+1.$ By Lemma 8, there exists a connected $k$-regular graph $C^{(3)}$ of order $n-(11q+10).$  We denote by $P^{(3)}$ the circular join of the eight graphs
$H_1^{(3)},\ldots, H_8^{(3)}$ where $H_1^{(3)}=H_2^{(3)}=K_{q+1},$ $H_3^{(3)}=H_8^{(3)}=K_q-PM,$ $H_4^{(3)}=H_7^{(3)}=K_{q+2}$ and $H_5^{(3)}=H_6^{(3)}=K_q.$
Denote $A^{(3)}=K_{3q+4}-HC.$ Let $V(A^{(3)})=\{u_1,u_2,\ldots, u_{3q+4}\}$ and $E(A^{(3)})=\{u_iu_j|\, j\ne i+1\}$ where $u_{3q+5}=u_1.$
Partition the vertex set of $A^{(3)}$ into eight subsets (four pairs)
\begin{align*}
&V_1^{(3)}=\{u_1\},\quad V_5^{(3)}=\{u_2\}\\
&V_2^{(3)}=\{u_{3+2j}|\,j=0,1,\ldots,(q-2)/2\},\quad V_6^{(3)}=\{u_{4+2j}|\,j=0,1,\ldots,(q-2)/2\}\\
&V_3^{(3)}=\{u_{q+1+2j}|\,j=1,\ldots,q/2\},\quad V_7^{(3)}=\{u_{q+2+2j}|\,j=1,\ldots,q/2\}\\
&V_4^{(3)}=\{u_{2q+1+2j}|\,j=1,\ldots,(q+2)/2\},\quad V_8^{(3)}=\{u_{2q+2+2j}|\,j=1,\ldots,(q+2)/2\}
\end{align*}
such that $A^{(3)}[V^{(3)}_j]=K_1$ for $j=1,5,$ $A^{(3)}[V^{(3)}_j]=K_{q/2}$ for $j=2,3,6,7,$ $A^{(3)}[V^{(3)}_j]=K_{(q+2)/2}$ for $j=4,8$ and  $A^{(3)}[V^{(3)}_1\cup V^{(3)}_5]=\overline{K_2},$ $A^{(3)}[V^{(3)}_s\cup V^{(3)}_{s+4}]=K_q-HP$ for $s=2,3,$ $A^{(3)}[V^{(3)}_4\cup V^{(3)}_8]=K_{q+2}-HP.$
Let $R^{(3)}_i=A^{(3)}[V^{(3)}_i],$ $i=1,\ldots, 8.$ Finally let $M_3$ be the graph obtained from the disjoint union $C^{(3)}+A^{(3)}+P^{(3)}$ by first
taking a nice connection of $R^{(3)}_i$ and $H^{(3)}_i$ for $i=1,\ldots, 8$ and then adding edges joining every vertex of $C^{(3)}$ to every vertex of
$A^{(3)}.$

Note that every vertex in $H^{(3)}_i$ has a unique neighbor in $R^{(3)}_i$ for $i=1,\ldots,8.$ It is easy to verify that ${\rm rad}(M_3)=2$ and
${\rm diam}(M_3)=4,$ the three graphs $C^{(3)},$ $A^{(3)}$ and $P^{(3)}$ are connected and $k$-regular, and they are the central subgraph, annular subgraph and
peripheral subgraph of $M_3$ which has order $n.$

Case (4). $k\equiv 1\,\, {\rm mod}\,\,3$ and $q$ is odd. Now $k=3q+1$ is even. Let $n\ge 14q+11.$  We have $n-(11q+9)\ge k+1.$ By Lemma 8, there exists a connected $k$-regular graph $C^{(4)}$ of order $n-(11q+9).$ We denote by $P^{(4)}$ the circular join of the eight graphs $H_1^{(4)},\ldots, H_8^{(4)}$ where
$H_1^{(4)}=H_5^{(4)}=K_q,$ $H_2^{(4)}=H_4^{(4)}=H_6^{(4)}=H_8^{(4)}=K_{q+1},$ $H_3^{(4)}=H_7^{(4)}=K_{q+1}-PM.$
Denote $A^{(4)}=K_{3q+3}-PM.$ We distinguish two subcases.

Subcase (4.1). $q=1.$  In this case $A^{(4)}=K_6-PM.$ Let $V(A^{(4)})=\{u_1,u_2,u_3,u_4,u_5,u_6\}$ and
$E(A^{(4)})=\{u_iu_j|\,i<j,\,\, (i,j)\ne (1,4),(2,5),(3,6)\}.$ Let $M_4$ be the graph obtained from the disjoint union $C^{(4)}+A^{(4)}+P^{(4)}$ by first
adding edges joining $u_i$ to each vertex in $H_j^{(4)}$ for $(i,j)=(1,1),\,(2,2),\,(2,3),\,(3,4),\,(4,5),\,(5,6),\,(5,7),(6,8)$ and
then adding edges joining every vertex of $C^{(4)}$ to every vertex of $A^{(4)}.$ It is easy to verify that ${\rm rad}(M_4)=2$ and
${\rm diam}(M_4)=4,$ the three graphs $C^{(4)},$ $A^{(4)}$ and $P^{(4)}$ are connected and $k$-regular, and they are the central subgraph, annular subgraph and
peripheral subgraph of $M_4$ which has order $n.$

Subcase (4.2). $q\ge 3.$  Partition the vertex set of $A^{(4)}$ into eight subsets $V_1^{(4)},\ldots, V_8^{(4)}$ such that $A^{(4)}[V^{(4)}_j]=K_1$ for $j=1,5,$
$A^{(4)}[V^{(4)}_j]=K_{(q-1)/2}$ for $j=2,6,$ $A^{(4)}[V^{(4)}_j]=K_{(q+1)/2}$ for $j=3,4,7,8$ and  $A^{(4)}[V^{(4)}_1\cup V^{(4)}_5]=\overline{K_2},$
$A^{(4)}[V^{(4)}_2\cup V^{(4)}_6]=K_{q-1}-PM,$ $A^{(4)}[V^{(4)}_s\cup V^{(4)}_{s+4}]=K_{q+1}-PM$ for $s=3,4.$ Let
$R^{(4)}_i=A^{(4)}[V^{(4)}_i],$ $i=1,\ldots, 8.$ Finally let $M_4$ be the graph obtained from the disjoint union $C^{(4)}+A^{(4)}+P^{(4)}$ by first
taking a nice connection of $R^{(4)}_i$ and $H^{(4)}_i$ for $i=1,\ldots, 8$ and then adding edges joining every vertex of $C^{(4)}$ to every vertex of
$A^{(4)}.$

Note that every vertex in $H^{(4)}_i$ has a unique neighbor in $R^{(4)}_i$ for $i=1,\ldots,8.$ It is easy to verify that ${\rm rad}(M_4)=2$ and
${\rm diam}(M_4)=4,$ the three graphs $C^{(4)},$ $A^{(4)}$ and $P^{(4)}$ are connected and $k$-regular, and they are the central subgraph, annular subgraph and
peripheral subgraph of $M_4$ which has order $n.$

Case (5). $k\equiv 2\,\, {\rm mod}\,\,3$ and $q$ is even.  Now $k=3q+2$ is even. Let $n\ge 14q+15.$  We have $n-(11q+12)\ge k+1.$ By Lemma 8, there exists a connected $k$-regular graph $C^{(5)}$ of order $n-(11q+12).$ We denote by $P^{(5)}$ the circular join of the eight graphs $H_1^{(5)},\ldots, H_8^{(5)}$ where
each $H_j^{(5)}=K_{q+1}$ for $j=1,\ldots,8.$ Denote $A^{(5)}=K_{3q+4}-PM.$ We distinguish two subcases.

Subcase (5.1). $q=0.$  In this case $k=2$ and $A^{(5)}=C_4.$ Let $V(A^{(5)})=\{u_1,u_2,u_3,u_4\}$ and $E(A^{(5)})=\{u_iu_{i+1}|\,i=1,2,3,4\}$ where $u_5=u_1.$
Let $M_5$ be the graph obtained from the disjoint union $C^{(5)}+A^{(5)}+P^{(5)}$ by first
adding edges joining $u_i$ to each vertex in $H_j^{(5)}$ for $(i,j)=(1,1),\,(2,2),\,(2,3),\,(2,4),\,(3,5),\,(4,6),\,(4,7),(4,8)$ and
then adding edges joining every vertex of $C^{(5)}$ to every vertex of $A^{(5)}.$ It is easy to verify that ${\rm rad}(M_5)=2$ and
${\rm diam}(M_5)=4,$ the three graphs $C^{(5)},$ $A^{(5)}$ and $P^{(5)}$ are connected and $k$-regular, and they are the central subgraph, annular subgraph and
peripheral subgraph of $M_5$ which has order $n.$

Subcase (5.2). $q\ge 2.$  Partition the vertex set of $A^{(5)}$ into eight subsets $V_1^{(5)},\ldots, V_8^{(5)}$ such that
$A^{(5)}[V^{(5)}_j]=K_1$ for $j=1,5,$ $A^{(5)}[V^{(5)}_j]=K_{q/2}$ for $j=2,3,6,7,$ $A^{(5)}[V^{(5)}_j]=K_{(q+2)/2}$ for $j=4,8$ and  $A^{(5)}[V^{(5)}_1\cup V^{(5)}_5]=\overline{K_2},$ $A^{(5)}[V^{(5)}_s\cup V^{(5)}_{s+4}]=K_q-PM$ for $s=2,3,$ $A^{(5)}[V^{(5)}_4\cup V^{(5)}_8]=K_{q+2}-PM.$  Let
$R^{(5)}_i=A^{(5)}[V^{(5)}_i],$ $i=1,\ldots, 8.$ Finally let $M_5$ be the graph obtained from the disjoint union $C^{(5)}+A^{(5)}+P^{(5)}$ by first
taking a nice connection of $R^{(5)}_i$ and $H^{(5)}_i$ for $i=1,\ldots, 8$ and then adding edges joining every vertex of $C^{(5)}$ to every vertex of
$A^{(5)}.$

Note that every vertex in $H^{(5)}_i$ has a unique neighbor in $R^{(5)}_i$ for $i=1,\ldots,8.$ It is easy to verify that ${\rm rad}(M_5)=2$ and
${\rm diam}(M_5)=4,$ the three graphs $C^{(5)},$ $A^{(5)}$ and $P^{(5)}$ are connected and $k$-regular, and they are the central subgraph, annular subgraph and
peripheral subgraph of $M_5$ which has order $n.$

Case (6). $k\equiv 2\,\, {\rm mod}\,\,3$ and $q$ is odd. Now $k=3q+2$ is odd. Let $n$ be even and $n\ge 14q+16.$ We have $n-(11q+13)\ge k+1.$ By Lemma 8, there exists a connected $k$-regular graph $C^{(6)}$ of order $n-(11q+13).$  We denote by $P^{(6)}$ the circular join of the eight graphs
$H_1^{(6)},\ldots, H_8^{(6)}$ where each $H_j^{(6)}=K_{q+1}$ for $j=1,\ldots,8.$ Denote $A^{(6)}=K_{3q+5}-HC.$ Let $V(A^{(6)})=\{u_1,u_2,\ldots, u_{3q+5}\}$ and $E(A^{(6)})=\{u_iu_j|\, j\ne i+1\}$ where $u_{3q+6}=u_1.$ Partition the vertex set of $A^{(6)}$ into eight subsets (four pairs)
\begin{align*}
&V_1^{(6)}=\{u_1\},\quad V_5^{(6)}=\{u_2\}\\
&V_2^{(6)}=\{u_{3+2j}|\,j=0,1,\ldots,(q-1)/2\},\quad V_6^{(6)}=\{u_{4+2j}|\,j=0,1,\ldots,(q-1)/2\}\\
&V_3^{(6)}=\{u_{q+2j}|\,j=2,\ldots,(q+3)/2\},\quad V_7^{(6)}=\{u_{q+1+2j}|\,j=2,\ldots,(q+3)/2\}\\
&V_4^{(6)}=\{u_{2q+1+2j}|\,j=2,\ldots,(q+3)/2\},\quad V_8^{(6)}=\{u_{2q+2+2j}|\,j=2,\ldots,(q+3)/2\}
\end{align*}
such that $A^{(6)}[V^{(6)}_j]=K_1$ for $j=1,5,$ $A^{(6)}[V^{(6)}_j]=K_{(q+1)/2}$ for $j=2,3,4,6,7,8$ and  $A^{(6)}[V^{(6)}_1\cup V^{(6)}_5]=\overline{K_2},$ $A^{(6)}[V^{(6)}_s\cup V^{(6)}_{s+4}]=K_{q+1}-HP$ for $s=2,3,4.$  Let $R^{(6)}_i=A^{(6)}[V^{(6)}_i],$ $i=1,\ldots, 8.$ Finally let $M_6$ be the graph obtained from the disjoint union $C^{(6)}+A^{(6)}+P^{(6)}$ by first taking a nice connection of $R^{(6)}_i$ and $H^{(6)}_i$ for $i=1,\ldots, 8$ and then adding edges joining every vertex of $C^{(6)}$ to every vertex of $A^{(6)}.$

Note that every vertex in $H^{(6)}_i$ has a unique neighbor in $R^{(6)}_i$ for $i=1,\ldots,8.$ It is easy to verify that ${\rm rad}(M_6)=2$ and
${\rm diam}(M_6)=4,$ the three graphs $C^{(6)},$ $A^{(6)}$ and $P^{(6)}$ are connected and $k$-regular, and they are the central subgraph, annular subgraph and
peripheral subgraph of $M_6$ which has order $n.$ This completes the proof. $\Box$

{\bf Remark 2.} In Theorem 9, the condition of being connected on metric subgraphs is essential. For example, Theorem 9 asserts that the smallest order
of a connected graph whose metric subgraphs are all cubic is $22.$ Let $Q$ be the graph obtained from the disconnected graph in Figure 2 by adding all the
edges $x_iy_j$ for $i=1,\ldots, 4$ and $j=1,\ldots, 6.$ Then $Q$ is a graph of order $18$ whose metric subgraphs are all cubic. The peripheral subgraph
of $Q$ is disconnected.
\vskip 3mm
\par
 \centerline{\includegraphics[width=3.7in]{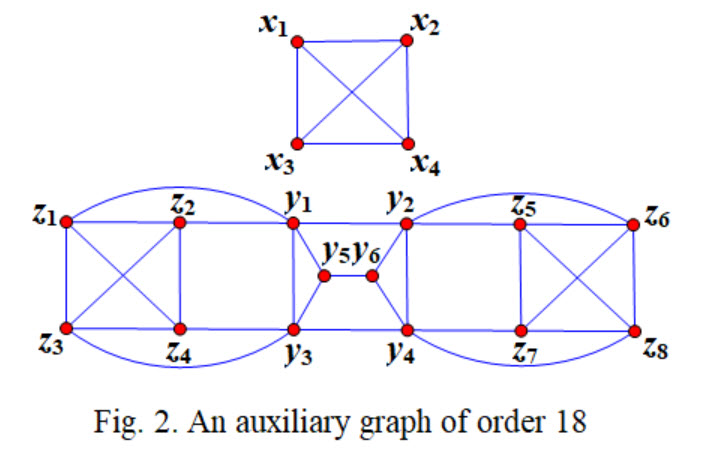}}
\par

Now we determine the smallest graphs whose metric subgraphs are all cycles.

{\bf Theorem 10.} {\it The minimum order of a connected graph whose metric subgraphs are all cycles is $15$ and there are exactly three such graphs
of order $15,$ all of which have size $35$ and are depicted in Figure 3. }

{\bf Proof.} The case $k=2$ of Theorem 9 asserts that the minimum order of a connected graph whose metric subgraphs are all cycles is $15.$
Let $G$ be a connected graph of order $15$ whose metric subgraphs are all cycles, and let $C,$ $A$ and $P$ be the center, annulus and periphery of $G$
respectively. Then $|C|+|A|+|P|=15.$ By the first four paragraphs of the proof of Theorem 9 we have $|C|\ge 3,$ $|A|\ge 4$ and $|P|\ge 8.$ Hence
$|C|=3,$ $|A|=4$ and $|P|=8.$ Denote by $H$ the peripheral subgraph of $G.$  Then $H=C_8$ and ${\rm rad}(H)=4.$

By Lemma 2, ${\rm rad}(G)\ge 2,$ ${\rm diam}(G)\ge 4$ and by Lemma 3, we have $4={\rm rad}(H)\ge {\rm diam}(G).$ Thus ${\rm diam}(G)=4$ and consequently
${\rm rad}(G)=2,$ since $A$ is nonempty. Combining the fact that the eccentricities of two adjacent vertices differ by at most $1$ and the condition
that ${\rm rad}(G)=2,$  we deduce that for any vertex $x\in C$ and any vertex $y\in P,$  $d(x,y)=2.$ Let $x,z,y$ be a path. Then $z\in A.$
Thus every vertex  in $P$ has a neighbor in $A.$ Note that the three metric subgraphs of $G$ are the cycles $C_3,$ $C_4$ and $C_8.$ Let $y^{\prime}$ be the antipodal vertex of $y$ on the even cycle $G[P].$ Then clearly $y^{\prime}$ is the unique eccentric vertex of $y$ in $G.$ The condition $d(y,y^{\prime})=4$ implies that
$y^{\prime}$ has a unique neighbor $z^{\prime}$ in $A$ which is the antipodal vertex of $z$ on the cycle $G[A].$ Since $y$ is the unique eccentric vertex of
$y^{\prime},$ $z$ is the unique neighbor of $y$ in $A.$ This shows that any vertex in $P$ has a unique neighbor in $A,$ implying that $|[P,\,A]|=8.$

We assert that every vertex in $A$ has a neighbor in $P.$ To see this, choose any vertex  $v\in A.$ Let $v^{\prime}$ be the antipodal vertex of $v$ in the even cycle $G[A].$ If $v$ has no neighbor in $P,$ then we would have $N(P)\subseteq A\setminus \{v\}$ and hence $e(v^{\prime})\le 2,$ a contradiction. Next we assert that every vertex in $C$ is adjacent to every vertex in $A.$ To show this, given any $u\in C$ and $v\in A$ we let $w$ be a neighbor of $v$ in $P.$
If $u$ and $v$ are nonadjacent, then $d(u,w)\ge 3,$ contradicting $e(u)=2.$ Now $|[C,\,A]|=12.$ Also, the three metric subgraphs $C_3,$ $C_4$ and $C_8$ contain
$3+4+8=15$ edges. Altogether $G$ has size $12+8+15=35.$

Combining the properties of $G$ deduced above we conclude that $G$ must be one of the three graphs depicted in Figure 3.
\vskip 3mm
\par
 \centerline{\includegraphics[width=5.7in]{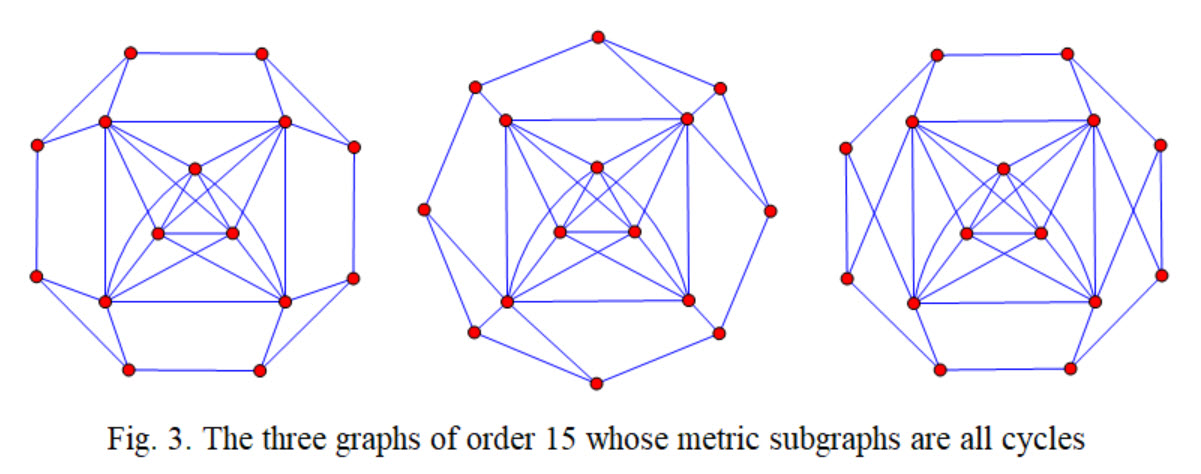}}
\par
Conversely, the metric subgraphs of each of these three graphs are all cycles. This completes the proof. $\Box$

{\bf Theorem 11.} {\it There exists a graph whose metric subgraphs are all isomorphic to a graph $H$ if and only if either $H$ is disconnected
or $H$ is connected and the radius of $H$ is at least $4.$ There exists a graph of order $n$ whose metric subgraphs are connected and pairwise isomorphic if and only if $n\ge 24$ and $n$ is divisible by $3.$ }

{\bf Proof.} Suppose $G$ is a graph whose metric subgraphs are all isomorphic to a graph $H.$ If $H$ is connected, then by Lemma 2, ${\rm diam}(G)\ge 4$ and by Lemma 3, ${\rm rad}(H)\ge {\rm diam}(G).$ Hence ${\rm rad}(H)\ge 4.$

Conversely, let $H$ be a given graph that is either disconnected or connected and ${\rm rad}(H)\ge 4.$ Recall that the eccentricity of any vertex in a disconnected graph is infinity. We will construct a graph $Q$ whose  metric subgraphs are all isomorphic to $H.$ Let $H_c,$ $H_a$ and $H_p$ be three
pairwise vertex-disjoint copies of $H,$ let $V(H_a)=\{x_1,\ldots, x_k\}$ and let $f$ be an isomorphism from $H_a$ to $H_p.$  Then the graph $Q$ is obtained
from $H_c +H_a +H_p$ by first taking the join of $H_c$ and $H_a$ and then adding the edges $x_if(x_i)$ for $i=1,\ldots,k.$ It is easy to verify that
$Q$ has radius $2$ and diameter $4,$ and the central subgraph, annular subgraph and peripheral subgraph of $Q$ are $H_c,$ $H_a$ and $H_p$ respectively.

Next we prove the second assertion of Theorem 11. Suppose $R$ is a graph of order $n$ whose metric subgraphs are connected and pairwise isomorphic. Let $W$ be the peripheral subgraph of $R.$ By Lemma 2, ${\rm diam}(R)\ge 4$ and by Lemma 3, $|W|\ge 2\, {\rm diam}(R).$ Hence $|W|\ge 2\times 4=8.$ It follows that $n=3 |W|\ge 24$ and $n$ is divisible by $3.$ Conversely, suppose $n\ge 24$ is an integer that is divisible by $3.$ Let $H$ be the path $P_{n/3}.$ By the first part
of Theorem 11 proved above, there exists a graph whose metric subgraphs are all isomorphic to $H.$ This completes the proof. $\Box$

Lemma 5 shows that connectedness conditions on the metric subgraphs of a graph yields some restriction on the order of the graph. Finally we pose the following question.

{\bf Question.} Let $k\ge 2$ be an integer. What is the smallest order of a graph whose metric subgraphs are all $k$-connected?

Using the information obtained in this paper, it is easy to show that the answer is $15$ when $k=2.$

\vskip 5mm
{\bf Acknowledgement.} This research  was supported by the NSFC grants 11671148 and 11771148 and Science and Technology Commission of Shanghai Municipality (STCSM) grant 18dz2271000.


\begin{thebibliography}{99}
\bibitem{1} H. Bielak and M.M. Syslo, Peripheral vertices in graphs, Studia Sci. Math. Hungar., 18(1983), no.2-4, 269-275.
\bibitem{2} F. Buckley, The central ratio of a graph, Discrete Math., 38(1982), 17-21.
\bibitem{3} F. Buckley and F. Harary, Distance in Graphs, Addison-Wesley Publishing Company, 1990.
\bibitem{4} G. Chartrand, D. Erwin, G. Johns and P. Zhang, Boundary vertices in graphs, Discrete Math., 263(2003), no.1-3, 25-34.
\bibitem{5} G. Chartrand, G. Johns, S. Tian and S.J. Winters, The interior and annulus of a graph, Congr. Numer., 102(1994), 57-62.
\bibitem{6} G. Chartrand, L. Lesniak and P. Zhang, Graphs and Digraphs, Sixth Edition, CRC Press, Boca Raton, 2016.
\bibitem{7} F. Harary and R.Z. Norman, The dissimilarity characteristic of Husimi trees, Ann. of Math., 58(1953), no.1, 134-141.
\bibitem{8} Y. Hu and X. Zhan, Possible cardinalities of the center of a graph, Bull. Malays. Math. Sci. Soc., 44(2021), 3629-3636.
\bibitem{9} B.M. Kim, Y. Rho, B.C. Song and W. Hwang, The maximum radius of graphs with given order and minimum degree, Discrete Math., 312(2012), 207-212.
\bibitem{10} L. Lesniak, Eccentric sequences in graphs, Period. Math. Hungar., 6(1975), no.4, 287-293.
\bibitem{11} D.B. West, Introduction to Graph Theory, Prentice Hall, Inc., 1996.
\end{thebibliography}
\end{document}